\documentclass[10pt,a4paper]{article}
\usepackage{amsthm,amsmath,amssymb,amscd,amsxtra,multicol}

\def\CC{\mathbb{C}}

\def\NN{\mathbb{N}}
\def\KK{\mathbb{K}}

\def\noi{\noindent}
\begin{document}
\newtheorem{defi}{Definition}
\newtheorem{exa}{Example}
\newtheorem{lem}{Lemma}
\newtheorem{pro}{Proposition}
\newtheorem{nota}{Notation}
\newtheorem{theo}{Theorem}
\newtheorem{coro}{Corollary}
\newtheorem{remark}{Remark}
\newtheorem{conjecture}{Conjecture}

{\bigbreak}
                          
\pagestyle{myheadings} \markboth{J. OJEDA}{\footnotesize{$p$-adic  meromorphic functions $f'P(f),\ g'P'(g)$ sharing a small function}}
\title{$p$-adic  meromorphic functions $f'P'(f),\ g'P'(g)$ sharing a small function}
\author { Kamal Boussaf, Alain Escassut and Jacqueline Ojeda}

\maketitle \thispagestyle{empty}

\begin{abstract}

Let $\KK$ be a complete algebraically closed p-adic field of characteristic zero. Let  $f,\ g$ be two   transcendental   meromorphic functions in the whole field $\KK$ or meromorphic functions in an open disk that are not  quotients of bounded analytic functions.
Let $P$ be a  polynomial of uniqueness for meromorphic functions in $\KK$ or in an open disk
  and let $\alpha$ be  a small meromorphic function  with regards to $f$ and $g$.
  If $f'P'(f)$ and $g'P'(g)$  share $\alpha$ counting multiplicity,  then we show that $f=g$ provided that the  multiplicity order of zeroes of $P'$ satisfy certain inequalities. If $\alpha$ is a Moebius function or a non-zero constant, we can obtain more general results on $P$. 

\end{abstract}

\noindent\\
\noindent
\noindent\footnotetext[1]{{\footnotesize {\sl Partially supported by CONICYT $N^\circ 9090014$ (Inserci\'on de Capital
Humano a la Academia)"}}}

 \noindent{\bf Introduction and Main Results}

\bigskip

Let  $\KK$ be an algebraically closed field of characteristic zero, complete for an ultrametric absolute value denoted by $|\ .\ |$. We denote by ${\cal A}({\KK})$ the $\KK$-algebra of entire functions in ${\KK}$, by ${\cal M}({\KK})$ the field of meromorphic functions in ${\KK}$, i.e. the field of fractions of ${\cal A}({\KK})$ and by ${\KK}(x)$ the field of rational functions.

Let $a \in {\KK}$ and $R \in ]0,+\infty[$. We denote by $d(a,R) $ the closed disk $\{ x \in {\KK} : |x - a| \leq  R \}$ and by  $d(a,R^-)$ the   open'' disk  $\{ x \in {\KK} : |x - a| < R \}$. We denote by ${\cal A}(d(a,R^-))$ the set of analytic functions in $d(a,R^-)$, i.e. the $\KK$-algebra of power series $\displaystyle {\sum_{n=0}^{\infty}}{a_n}(x - a)^n$ converging in $d(a,R^-)$ and by ${\cal M} (d(a,R^-))$ the field of meromorphic functions inside $d(a,R^-)$, i.e. the field of fractions of  ${\cal A} (d(a,R^-))$. Moreover, we denote by ${\cal A}_b (d(a,R^-))$ the ${\KK}$ - subalgebra of ${\cal A} (d(a,R^-))$ consisting of the bounded analytic functions in $d(a,R^-)$, i.e. which satisfy $\displaystyle \sup_{n \in \NN} |a_n|R^n < + \infty$ . And we denote  by ${\cal M}_b (d(a,R^-))$ the field of fractions of ${\cal A}_b (d(a,R^-))$. Finally, we denote by ${\cal A}_u (d(a,R^-))$ the set of unbounded analytic functions in $d(a,R^-)$, i.e. ${{\cal A}(d(a,R^-)) \setminus {\cal A}_b (d(a,R^-))}$. Similarly, we set ${\cal M}_u (d(a,R^-)) = {\cal M}(d(a,R^-)) \setminus {\cal M}_b (d(a,R^-))$. 

\medskip

The problem of value sharing a small function   by  functions of the form $f'P'(f)$ was examined first when  $P$ was just of the form $x^n$ [7],  [18], [24].  More recently it was examined when $P$ was a  polynomial such that $P'$ had exactly two distinct zeroes  [15],  [17], [20], both in complex analysis and in p-adic analysis. In [15], [17] the   functions where meromorphic on $\CC$, with a 
 small function that was a constant or the identity. In [20], the problem was considered for analytic functions in the field $\KK$: on one hand for entire functions and on the other hand for  unbounded analytic functions in an open disk.

 Here we consider functions $f,\ g\in {\cal M}(\KK)$ or $f,\ g\in {\cal M}(d(a,R^-))$ and ordinary polynomials $P$:   we must only assume certain hypotheses on the multiplicity order of the zeroes of $P'$.    
   The method for the  various   theorems  we will show is the following: assuming that $f'P'(f)$ and $g'P'(g)$ share a small function, we first prove that $f'P'(f)=g'P'(g)$. Next, we derive  $P(f)=P(g)$. And then, when $P$ is a polynomial of uniqueness for the functions we consider, we can conclude $f=g$.

\medskip

Now, in order to define small functions, we have to 
 briefly recall the definitions of the classical Nevanlinna theory in the field $\KK$ and a few specific properties of ultrametric analytic or meromorphic functions.
\smallskip

Let $\log$ be a real logarithm function of base $>1$ and let $f \in {\cal M}(\KK)$ \big(resp. $f \in {\cal M}(d(0,R^-))$\big) having no zero and no pole at 0. Let $r \in ]0,+\infty[$ \big(resp. $r \in ]0,R[$\big) and let $\gamma \in d(0,r)$. If $f$ has a zero of order $n$ at $\gamma$, we put $\omega_\gamma(h)=n$. If $f$ has a pole of order $n$ at $\gamma$, we put $\omega_\gamma(f)=-n$ and finally, if $f(\gamma) \not= 0, \infty$, we put $\omega_\gamma(f)=0$
\smallskip

We denote by $Z(r,f)$ the {\it counting function of zeroes of $f$} in $d(0,r)$, counting multiplicity, i.e. we set 
$$Z(r,f) = \sum_{\omega_\gamma(f)>0, \ |\gamma|\leq r} \omega_\gamma(f) (\log r - \log |\gamma|).$$

In the same way, we set $\displaystyle N(r,f) = Z\Big(r,\frac{1}{f}\Big)$   to denote the {\it counting function of poles of $f$} in $d(0,r)$, counting multiplicity.
\smallskip

For $f \in {\cal M}(d(0,R^-))$ having no zero and no pole at 0, the {\it Nevanlinna function} is defined by $\displaystyle T(r,f) = \max \big\{ Z(r,f) + \log |f(0)|, N(r,f) \big\}$.

\smallskip  

Now,  we must recall the definition of a {\it small function} with respect to a meromorphic function and some pertinent properties.

\medskip

\noindent {\bf Definition.}  
Let $f \in {\cal M}(\KK)$ \big(resp. let $f \in {\cal M}(d(0,R^-))$\big) such that $f(0) \not=0, \infty$. A function $\alpha \in {\cal M}(\KK)$  \big(resp. $\alpha \in {\cal M}(d(0,R^-))$\big) having no zero and no pole at 0 is called {\it a small function with respect to $f$}, if it satisfies $\displaystyle \lim_{r \to +\infty}\frac{T(r,\alpha)}{T(r,f)} = 0 \ \ \Big({\rm resp.} \ \lim_{r \to R^-}\frac{T(r,\alpha)}{T(r,f)} = 0 \Big)$.

If 0 is a zero or a pole of $f$ or $\alpha$, we can make a change of variable such that the new origin is not a zero or a pole for both $f$ and $\alpha$. Thus it is easily seen that the last relation does not really depend on the origin.

\smallskip

We denote by ${\cal M}_f(\KK)$ \big(resp. ${\cal M}_f(d(0,R^-))$\big) the set of small meromorphic functions with respect to $f$ in $\KK$ \big(resp. in $d(0,R^-)$\big).

\smallskip  
 Let us remember the following definition.
 
\medskip

  \noindent {\bf Definition.} 
Let $f,g, \alpha \in {\cal M}(\KK)$ \big(resp. let $f,g, \alpha \in {\cal M} (d(0,R^-))$\big). We say that $f$ and $g$ {\it share the function $\alpha$ C.M.}, if $f - \alpha$  and $g - \alpha$ have the same zeroes with the same multiplicity in $\KK$ \big(resp. in $d(0,R^-)$\big).

Recall that a polynomial $P\in \KK[x]$ is called {\it a  polynomial of uniqueness} for a class of functions $\cal F$ if for any two functions $f,\ g\in {\cal F}$  the property $P(f)=P(g)$ implies $f=g$. 

 \bigskip
 
Actually, in   a  p-adic field, we can obtain  various results, not only  for functions defined in the whole field $\KK$ but also   for functions defined inside an open disk because the p-adic Nevanlinna Theory works inside a disk, for functions of ${\cal M}_u(d(0,R^-))$.

 \bigskip

We can now state our main theorems on the problem $f'P'(f), \ g'P'(g)$ share a  small function.

\medskip

\begin{theo}
Let $P$ be a polynomial of uniqueness for ${\cal M}(\KK)$,  let $ \displaystyle{P'=b(x-a_1)^n\prod
_{i=2}^l(x-a_i)^{k_i}}$ with $b\in \KK^*$,   $l\geq 2$,  $k_i\geq k_{i+1}, \ 2\leq i\leq l-1$ and let $k=\sum_{i=2}^lk_i$.  Suppose $P$ satisfies the following conditions:

$n\geq 10 +\displaystyle {\sum_{i=3}^{l}\max(0, 4-k_i)+\max(0, 5-k_2)}, $

 $n\geq k+2,$ 

if $l=2,$  then   $  n  \neq  2k,\ 2k+1,\ 3k+1$, 

if $l=3, $  then $  n\neq 2k+1,\ 3k_i-k\ \forall  i=2,3$.

Let $f, g \in {\cal M}(\KK)$ be transcendental and let $\alpha \in {\cal M}_f(\KK)\cap {\cal M}_g(\KK)$ be non-identically zero.
  If $f'P'(f)$ and $g'P'(g)$ share $\alpha$ C.M., then $f=g$.

 \end{theo}

\begin{theo}
Let $P$ be a polynomial of uniqueness for ${\cal M}(\KK)$,  let $ \displaystyle{P'=b(x-a_1)^n\prod
_{i=2}^l(x-a_i)^{k_i}}$ with $b\in \KK^*$,  $l \geq 2$,   $k_i\geq k_{i+1}, \ 2\leq i\leq l-1$ and let $k=\sum
_{i=2}^lk_i$.
 Suppose $P$ satisfies the following conditions:

$n\geq 9 +\displaystyle {\sum
_{i=3}^{l}\max(0, 4-k_i)+\max(0, 5-k_2)}, $

 $n\geq k+2,$ 

if $l=2,$ then $ n  \neq  2k,\ 2k+1,\ 3k+1$, 

if $l=3,  $ then $  n\neq 2k+1,\  3k_i-k\ \forall  i=2,3$.

Let $f, g \in {\cal M}(\KK)$ be transcendental and let $\alpha $ be a Moebius function. If $f'P'(f)$ and $g'P'(g)$ share $\alpha$ C.M.,  then $f=g$.

 \end{theo}

\begin{theo}
Let $P$ be a polynomial of uniqueness for ${\cal M}(\KK)$,  let $ \displaystyle{P'=b(x-a_1)^n\prod
_{i=2}^l(x-a_i)^{k_i}}$ with $b\in \KK^*$,  $l \geq 2$,  $k_i\geq k_{i+1}, \ 2\leq i\leq l-1$ and let $k=\sum
_{i=2}^lk_i$.
 Suppose $P$ satisfies the following conditions:

 $n\geq k+2$,
 
$n\geq 9 +\displaystyle {\sum
_{i=3}^{l}\max(0, 4-k_i)+\max(0, 5-k_2)}. $

Let $f, g \in {\cal M}(\KK)$ be transcendental and let $\alpha  $ be   a non-zero  constant.  If $f'P'(f)$ and $g'P'(g)$ share $\alpha$ C.M.,  then $f=g$.

 \end{theo}

 \begin{theo} Let $a\in K $ and $R>0$.  Let $P$ be a polynomial of uniqueness for ${\cal M}_u(d(a,R^-))$ and  let $ \displaystyle{P'=b(x-a_1)^n\prod
_{i=2}^l(x-a_i)^{k_i}}$ with $b\in \KK^*$,  $l \geq 2,  \ k_i\geq k_{i+1}, \ \ 2  \leq  i\leq l-1$  and let $k=\sum
_{i=2}^lk_i$.   Suppose $P$ satisfies the following conditions:

$n\geq 10 +\displaystyle {\sum
_{i=3}^{l}\max(0, 4-k_i)+\max(0, 5-k_2)}, $

 $n\geq k+3,$ 

if $l=2, $ then $   n  \neq  2k,\ 2k+1,\ 3k+1, $

if $l=3,  $ then $ n\neq 2k+1, \  3k_i-k\ \forall  i=2,3$.

Let $f, g \in {\cal M}_u(d(a,R^-))$  and let $\alpha \in {\cal M}_f(d(a,R^-))\cap {\cal M}_g(d(a,R^-))$ be non-identically zero.   If $f'P'(f)$ and $g'P'(g)$ share $\alpha$ C.M.,  then $f=g$.

 \end{theo}

\begin{theo}
Let $P$ be a polynomial of uniqueness for ${\cal M}(\KK)$  such that $P'$ is of the form

\noindent $\displaystyle{b(x-a_1)^n\prod
_{i=2}^l(x-a_i)}$ with   $l\geq 3$ ,  $b\in \KK^*$, satisfying:

 $n\geq l+10$, 

if $l=3$, then $n\neq 2l-1$. 

\noindent  Let $f, g \in {\cal M}(\KK)$  be transcendental   and let $\alpha \in {\cal M}_f(\KK)\cap {\cal M}_g(\KK)$ be non-identically zero.  If $f'P'(f)$ and $g'P'(g)$ share $\alpha$ C.M.,  then $f=g$.

 \end{theo}

\begin{theo} Let $a\in K $ and $R>0$.
Let $P$ be a polynomial of uniqueness for ${\cal M}_u(d(a,R^-))$  such that $P'$ is of the form 
 $\displaystyle{P'=b(x-a_1)^n\prod
_{i=2}^l(x-a_i)}$ with $l\geq 3$,  $b\in \KK^*$ satisfying:

   $n\geq l+10$,
   
     if $l=3$, then $n\neq 2l-1$.

\noindent  Let $f, g \in {\cal M}_u(d(a,R^-))$    and let $\alpha \in {\cal M}_f(d(a,R^-))\cap {\cal M}_g(d(a,R^-))$ be non-identically zero.  If $f'P'(f)$ and $g'P'(g)$ share $\alpha$ C.M.,  then $f=g$.
 
  \end{theo}

\begin{theo}  
Let $P$ be a polynomial of uniqueness for ${\cal M}(\KK)$  such that $P'$ is of the form 

\noindent $\displaystyle{P'=b(x-a_1)^n\prod
_{i=2}^l(x-a_i)}$ with $l\geq 3$, $b\in \KK^*$ satisfying

  $n\geq l+9$,
  
  if $l=3$, then $n\neq 2l-1$. 

Let $f, g \in {\cal M}(\KK)$   be transcendental  and let $\alpha$ be   a  Moebius function.
  If $f'P'(f)$ and $g'P'(g)$ share $\alpha$ C.M.,  then $f=g$. 

 \end{theo}

\begin{theo}  
Let $P$ be a polynomial of uniqueness for ${\cal M}(\KK)$  such that $P'$ is of the form 

\noindent $\displaystyle{P'=b(x-a_1)^n\prod
_{i=2}^l(x-a_i)}$ with $l\geq 3$, $b\in \KK^*$ satisfying
$n\geq l+9$.

Let $f, g \in {\cal M}(\KK)$   be transcendental  and let $\alpha$ be  a non-zero constant. 
  If $f'P'(f)$ and $g'P'(g)$ share $\alpha$ C.M.,  then $f=g$. 

 \end{theo}

\begin{theo}
Let $f,\   g \in {\cal M}(\KK)$ be transcendental and let $\alpha \in {\cal M}_f(\KK)\cap {\cal M}_g(\KK)$ be non-identically zero. Let $a \in \KK \setminus \{0\}$. If  $f'f^n(f-a)$ and $g'g^n(g-a)$ share the function $\alpha$ C.M. and if $n  \geq 12$, then either $f=g$ or there exists $h\in{\cal M}( \KK)$  such that $\displaystyle f= \frac{a(n+2)}{n+1}\Big(\frac{h^{n+1}-1}{h^{n+2}-1}\Big)h$ and $\displaystyle g= \frac{a(n+2)} {n+1}\Big(\frac{h^{n+1}-1}{h^{n+2}-1}\Big)$. Moreover,
  if $\alpha$ is a  constant  or a  Moebius function, then the conclusion holds whenever $n\geq 11$.  

\end{theo}

Inside an open disk, we have a version similar to the general case in the whole field. 
   
\begin{theo}
Let $f, \  g \in {\cal M}_u(d(0,R^-))$, and let $\alpha \in {\cal M}_f(d(0,R^-))\cap {\cal M}_g(d(0,R^-))$  be non-identically zero. Let $a \in \KK \setminus \{0\}$. If  $f' f^n(f-a)$ and $g'g^n(g-a)$ share the function $\alpha$ C.M. and $n  \geq 12$, then either $f=g$ or there exists  $h \in {\cal M}(d(0,R^-))$ such that $\displaystyle f= \frac{a(n+2)} {n+1}\Big(\frac{h^{n+1}-1}{h^{n+2}-1}\Big)h$ and $\displaystyle g= \frac{a(n+2)} {n+1}\Big(\frac{h^{n+1}-1}{h^{n+2}-1}\Big)$.   

\end{theo}

\noi{\bf Specific Lemmas}

\begin{lem} Let 
 $f,\ g\in {\cal M}(\KK)$ be transcendental   (resp.  
$f,\ g\in {\cal M}_u(d(0,R^-))$). Let $P(x)=x^{n+1}Q(x)$ be a polynomial   such that  $n\geq deg(Q)+2$ (resp.   $n\geq deg(Q)+3$). If $P'(f)f'=P'(g)g'$ then $P(f)=P(g)$.

\end{lem}

\bigskip

\begin{lem} Let  $Q(x)=(x-a_1)^n\prod
_{i=2}^l(x-a_i)^{k_i}\in \KK[x]$ ($a_i\not=a_j, \  \forall i\not =j$) with $l\ge 2$ and  $n\geq \max\{k_2,..,k_l\}$ and  let $k=\sum
_{i=2}^lk_i$. Let 
 $f,\ g\in {\cal M}(\KK)$ be transcendental   (resp.  
$f,\ g\in {\cal M}_u(d(0,R^-))$) such that $\theta=Q(f)f'Q(g)g'$ is a small function with respect to $f$ and $g$.  We have the following :

If $l=2$ then $n$ belongs to $\{k,\ k+1,\ 2k,\ 2k+1,\ 3k+1\}$.

If $l=3$ then $n$ belongs to $\{\frac{k}{2},k+1,2k+1,3k_2-k,..,3k_l-k\}$.

If $l\geq 4$ then $n=k+1$.
 
If    $\theta$ is a constant and  $f,\ g\in {\cal M}(\KK)$  then $n=k+1$.
\end{lem}

\medskip \begingroup
Kaml BOUSSAF,  Alain ESCASSUT

Laboratoire de Mathematiques
UMR 6620

Universit\'e Blaise Pascal

Les C\'ezeaux

63171    AUBIERE

FRANCE

\bigskip

kamal.boussaf@math.univ-bpclermont.fr
\medskip

alain.escassut@math.univ-bpclermont.fr
\medskip

Jacqueline OJEDA

Departamento de Matematica

Facultad de Ciencias F'sicas y Matematicas

Universidad de Concepcion

\bigskip

mail: jacqojeda@udec.cl  
 
  Jacqueline.Ojeda@math.univ-bpclermont.fr
\endgroup

\end{document}